\documentclass[12pt]{article}

\usepackage[latin1]{inputenc}
\usepackage[T1]{fontenc}
\usepackage[french]{babel}
\usepackage{amsfonts}
\usepackage{euscript}

\newtheorem{conj}{Conjecture}[section]
\newtheorem{theo}[conj]{Theorem}
\newtheorem{prop}[conj]{Proposition}
\newtheorem{coro}[conj]{Corollary}
\newtheorem{lemm}{Lemma}[section]

\begin{document}

\newcommand{\as}{${\cal A}_6\ $}
\newcommand{\mt}{${\cal M}_3\ $}
\newcommand{\Pro}{{\bf I \hspace{-3pt} P}}
\newcommand{\p}{{\bf I \hspace{-3pt} P}}
\newcommand{\C}{\mathbb{C}}
\newcommand{\oo}{\mathbb{O}}
\newcommand{\n}{\mathfrak n}
\newcommand{\B}{\mathfrak b}
\newcommand{\nd}{$\frac{n}{2}\,$}
\newcommand{\nq}{$\frac{n}{4}\,$}
\newcommand{\f}{${\cal F}(Q_1)\ $}
\newcommand{\fp}{${\cal F}'(Q_1)\ $}
\newcommand{\fx}{${\cal F}_x(Q_1)\ $}
\newcommand{\fpx}{${\cal F}_x'(Q_1)\ $}
\newcommand{\se}{Sec(X)-X}
\newcommand{\sv}{Severi variety }
\newcommand{\vs}{Severi variety }
\newcommand{\g}{\mathfrak g}
\newcommand{\h}{\mathfrak h}

\newcommand{\dem}{\underline {\bf Proof :} }
\newcommand{\rem}{\underline {\bf Remark :} }
\newcommand{\fin}{\begin{flushright}  $\bullet$ \end{flushright}}
\newcommand{\para}{\vspace{1pt} \ \\}
\newcommand{\Para}{\vspace{15pt} \ \\}

\newcommand{\surmap}{\longrightarrow \hspace{-.5 cm} \longrightarrow}
\newcommand{\pisur}[2]{
\pi_1(#1) \surmap \pi_1(#2)   } 

\newcommand{\fonction}[5]{
\begin{array}{rrcll}
#1 & : & #2 & \rightarrow & #3 \\
   &   & #4 & \mapsto     & #5
\end{array}  }
\newcommand{\fonc}[3]{
#1 : #2 \rightarrow #3 }

\newcommand{\DC}[8]{
\begin{array}{ccc}
#1             & \stackrel{#2}{\longrightarrow} & #3 \\
\downarrow #4  &                                & \downarrow #5 \\
#6             & \stackrel{#7}{\longrightarrow} & #8
\end{array} }

\newcommand{\matdd}[4]{
\left (
\begin{array}{cc}
#1 & #2  \\
#3 & #4
\end{array}
\right )   }

\newcommand{\mattt}[9]{
\left (
\begin{array}{ccc}
  #1 & #2 & #3 \\
  #4 & #5 & #6 \\
  #7 & #8 & #9
\end{array}
\right )   }

\newcommand{\vectq}[4]{
\left (
\begin{array}{r}
#1 \\
#2 \\
#3 \\
#4
\end{array}
\right )   }

\title{Severi varieties}
\author{Pierre-Emmanuel Chaput}
\date{January 2001}
\maketitle

{\def\thefootnote{\relax}
\footnote{\hskip-0.6cm
{\it AMS mathematical classification \/}: 14M07,14M17,14E07.\newline
}}

\begin{center}
\bf{Abstract}
\end{center}

R. Hartshorne conjectured and F. Zak proved (cf \cite[p.9]{lazarsfeld})
that any smooth non-degenerate complex algebraic variety
$X^n \subset \p^m$ with $m<\frac{3}{2}n+2$ satisfies 
$Sec(X)=\p^m$ ($Sec(X)$ denotes the secant variety of $X$; when
$X$ is smooth it is simply the union of all the secant and tangent 
lines to $X$).
In this article, I deal with the limiting case of this theorem, namely 
the Severi varieties, defined by the conditions $m=\frac{3}{2}n+2$ and
$Sec(X)\not = \p^m$. I want to give a different proof of a theorem
of F. Zak classifying
all Severi varieties. F. Zak proves that there exists only four 
Severi varieties
and then realises a posteriori that all of them are homogeneous; here I will
work in another direction: I prove a priori that any Severi variety is
homogeneous and then deduce more quickly their classification, satisfying 
R. Lazarsfeld et A. Van de Ven's wish \cite[p.18]{lazarsfeld}.
By the way, I give a very brief proof of the fact that the derivatives of
the equation of $Sec(X)$, which is a cubic hypersurface,
determine a birational morphism of $\p^m$.

I wish to thank Laurent Manivel for helping me in writing this article.

\section{Preliminary facts about Severi varieties and an example}

In this section, I state results of F. Zak (cf \cite[p.19]{lazarsfeld}). The 
reader can find the details of their proofs in my preprint written in french
(cf \cite{prepub}).\\
Let $X$ be a Severi variety and $Sec(X)$ its secant variety.
Let, for any $P\in Sec(X)-X$,\\
$Q_P:=\{x\in X:(xP)$ secant or tangent to X$\}$ and
$\Sigma_P:=\bigcup_{x\in Q_P}(xP)$, 
the union of all
secant and tangent lines through P. It may also be described 
as the cone on $P$ with basis $Q_P$.

\begin{theo}{For any $P\in Sec(X)-X$,
\begin{itemize}{
\item{a: $Q_P$ is a smooth \nd-dimensional quadric, and $\Sigma_P$ is
a (\nd +1)-dimensional linear space.}
\item{b: $\Sigma_P \cap X = Q_P$}
\item{c: $\Sigma_P -X     = \{P'\in \se:T_{P'}Sec(X)=T_PSec(X)\}$  }
\item{d: $\forall P'\in \se,Q_P=Q_{P'} \Leftrightarrow P'\in\Sigma_P$}}
\end{itemize} \label{thQP} }
\end{theo}

\begin{theo}{$Sec(X)$ is a cubic hypersurface singular exactly on $X$.  
\label{thsec} }
\end{theo}

Their proofs imply the

\begin{lemm}{Let $P\in \se$ and $M$ a (\nd+2)-linear space
containing $\Sigma_P$ and not contained in $T_PSec(X)$. Then their exists
$x\in X-\Sigma_P$ such that
\begin{itemize}{
\item{(i)  $M\cap X=Q_P\cup\{x\}$}
\item{(ii) $M\cap Sec(X)=\Sigma_P \cup C(x,Q_P)$
($C(x,Q_P)$ is the cone on
$x$ with base $Q_P$)}  }
\end{itemize}    \label{M}  }
\end{lemm}

Now let's study the example of $\p^2 \times \p^2$, which is the variety of
$3\times 3$-matrices of rank 1 in $\p^8=\p[{\cal M}_3 (\C)]$. It is indeed a
\sv since matrices in $Sec(X)$ have rank at most 2; 
and we can see from this that $Sec(X)$ is defined by the determinant, 
an equation 
of degree 3 (cf theorem \ref{thsec}). If $P=M+N$, 
with $M,N\in X$ and $P\in \se$,
we have $rk(P)=2$, $Ker(M)\cap Ker(N)=Ker(P)$ 
and $Im(M)+Im(N)=Im(P)$. Let's consider the case where
$P=\matdd{Id}{0}{0}{0}$, we can then easily compute that:
$$\Sigma_P=\left \{ \mattt{a}{b}{0} {c}{d}{0} {0}{0}{0} \right \}$$ and
$$Q_P=\Sigma_P \cap \{ad-bc=0\}$$
\para

This illustrates theorem \ref{thQP}. We now go on analysing this particular
example so as to prepare the general method: let $F$ be the symmetric
trilinear form such that $F(M,M,M)=det(M)$, namely for 
$M_i=(C_{i,j})_{1\leq j \leq 3}$
any matrix considered as 3 column-vectors ($i=1,2,3$),
$$F(M_1,M_2,M_3)=\frac{1}{6}\sum_{\sigma \in S_3} \det(C_{\sigma(i),i})_i$$
Identifying a matrix with its orthogonal hyperplane for the canonical 
bilinear form, one can easily compute that
$\fonc{\tilde{G}}{M}{F(M,M,.)},M\in {\cal M}_3(\C)$ is
just the comatrix application and then the fonction $G$ defined by
$G(M)=\frac{F(M,M,.)}{F(M,M,M)}$ is an involution out of $Sec(X)$ as
it is the inverse application.
On the other side, we may define a regular fonction on $\p^m - X$ 
which coincide with $G$ on $\p^m - Sec(X)$ (by setting $G(M)=\mbox{Com}(M)$)
and this morphism verifies $G[Sec(X)]=X$. We are going to see that all of this
can be proved in the general case, except for the fact that we won't
identify a priori $\p^m$ and ${\p^m}^*$. 

We are also going to show that any \vs is homogeneous. For any variety $Z$,
let $Z^*$ be its dual variety, that is the variety of all hyperplanes tangent
to $Z$ at one point. In the example we have
considered, any invertible matrix yields a linear isomorphism $L_M$ between
$\p^8$ and ${\p^8}^*$ (namely $L_M(B)$ is the hyperplane in $\C^9$ orthogonal
to the matrix $-M^{-1}BM^{-1}$ for the canonical scalar product \footnote{
$L_M(B)=dG_M(B)$; in the second part, we are going to see the interest of 
this fonction.}). This isomorphism maps
$X$ on $Sec(X)^*$ and, for $M$ and $N$ varying, the endomorphisms of $\p^8$ 
${(L_N)}^{-1} \circ L_M$ restrict
to a family of endomorphisms acting transitively on $X$, proving
$X$'s homogeneity. The same will be shown to be true in the general case.

\section{Homogeneity of Severi varieties}

Let $X^n\subset \p^m$ be any \sv, 
$Y=Sec(X)^*\subset {\p^m}^*$, and as before
$F$ a symmetric trilinear form 
such that $v\in Sec(X) \Leftrightarrow F(v,v,v)=0$. We will simplify the
expression $F(v,v,v)$ in $F(v)$. From now on, all
elements will be considered as elements of $\C^{m+1}$ (and not 
$\p^m$) and we will denote the same way varieties in $\p^m$ and
their cones in $\C^{m+1}$. For 
$w_0\not \in Sec(X)$ and $w\in \C^{m+1}$ let us denote by 
$L_{w_0}(w)$ the linear form 
$2F(w_0)F(w_0,w,.)-3F(w_0,w_0,w)w_0^*$, where $w_0^*$ stands for the
linear form
$F(w_0,w_0,.)$. Let us also notice that $L_{w_0}(w)$ and the differential 
$D_{w_0}G(w)$ are colinear if $G$ stands for the rational map
$w\mapsto \frac{w^*}{F(w)}$ we have just mentionned.

Let $w_0$ be fixed and $x\in X$ such that $w_0^*(x)\not = 0$. I am going to
explain the geometric meaning of $L_{w_0}(x)$. Taking into account the fact
that $x^*=0$ ($X$
is the singular locus of $Sec(X)$), we know that
$x+\lambda w_0 \in Sec(X)$ if and only if
$\lambda=0 \mbox{ or } \lambda=\lambda_x:=
-\frac{w_0^*(x)}{F(w_0)} $.\\
On the other direction, let $H_{w_0}(x)$ be the tangent space to $Sec(X)$ at
$x+\lambda_x w_0$. The equation of this hyperplane is 
$F(x+\lambda_x w_0,x+\lambda_x w_0,.)=0$, or $L_{w_0}(x)(.)=0$.

\para
So $D_{w_0}G(x)$ is the equation of the tangent space to
$Sec(X)$ at the other intersection point of $(xw_0)$ with $Sec(X)$ if 
$x$ stands outside the hyperplane 
$w_0^*=0$. So $D_{w_0}G(X)\subset Y$. In the opposite direction, if $P\in \se$
verifies $w_0 \not \in H_P:=T_PSec(X)$, 
we know from lemma \ref{M}
that there exists a
point $x_P$ in $X\cap (\Sigma_P+w_0) - Q_P$, and then $x_P$ 
is such that $D_{w_0}G(x_P)=H_P$.
As $Y$
is known to be non degenerate, $D_{w_0}G(X)$ 
contains an open subset of $Y$ and we have the

\begin{prop}{If $w_0\not \in Sec(X)$, $D_{w_0}G$ is a linear isomorphism
from $\p^m$ to ${\p^m}^*$ such that $D_{w_0}G(X)=Y$.}
\end{prop}

\begin{coro}{$X$ is homogeneous.}
\end{coro}

\dem Let $x,x'\in X$; it is enough to find $w,w'\not \in Sec(X)$ such that
$D_wG(x)=D_{w'}G(x')$. Let $H_P$ containing neither $x$ nor $x'$; $w$
in $\Sigma_P + x$ and $w'$ in $\Sigma_P + x'$ achieve it.
\fin

\section{Classification of Severi varieties}

We now want to prove Zak's classification theorem (cf \cite{lazarsfeld})

\begin{theo}{There are only four Severi varieties, namely:
\begin{itemize}  {
\item{n=2 :  $X={\cal V}\subset\p^5$, the Veronese surface}
\item{n=4 :  $X=\p^2 \times \p^2 \subset \p^8$}
\item{n=8 :  $X=G(2,6) \subset \p^{14}$}
\item{n=16:  $X={\cal E}_6 \subset \p^{26}$}  }
\end{itemize}  \label{classification} }
\end{theo}

The reader may find in \cite[p.13]{lazarsfeld} elements to understand the
construction of ${\cal E}_6$ and the fact that it is homogeneous under $E_6$.
\para

\dem Let $V=\C^{m+1}$. Let $\cal H$ be the group of automorphisms of $\p V$
preserving $X$. As $\cal H$ acts transitively 
on the non-degenerate homogeneous projective variety $X$,
it is a semi-simple subgroup of $PGL(V)$. Let $\pi$ be the projection
$SL(V)\rightarrow PSL(V)$ and $\cal G$ the identity component of 
$\pi^{-1}({\cal H})$. $\cal G$
is again a semi-simple subgroup, and $V$ is an irreducible
$\cal G$-module.

\para
If $\cal G$ is not simple, we may write
${\cal G}={\cal G}_1 \times {\cal G}_2$
with non-trivial ${\cal G}_1$ and ${\cal G}_2$ acting on $V_1$ and $V_2$ 
such that
$V=V_1 \otimes V_2$. Let $n_i:=\dim X_i$ and $n_i+\delta_i:=\dim V_i$ 
(then $\delta_i \geq 1$).
For $X=X_1 \times X_2$ to be a \vs,
we need $\frac{3}{2}(n_1+n_2)+3=(n_1+\delta_1)\times(n_2+\delta_2)$,
which leads to:
$(n_1n_2) + (\delta_2-\frac{3}{2})n_1 + (\delta_1-\frac{3}{2})n_2 
+ (\delta_1 \delta_2 -3) = 0$.

If $\delta_1 \geq 2$ and $\delta_2 \geq 2$, all the terms of this sum are
positive; we can then suppose $\delta_2=1$, which gives
$(\delta_1-\frac{3}{2})n_2 + (n_2-\frac{1}{2})n_1 + \delta_1 - 3 =0$ and so
$\delta_1 \leq 2$.

If $\delta_1=2$, the equation gives 
$n_1(n_2-\frac{1}{2}) + \frac{n_2}{2} - 1=0$
so $n_1=n_2=1$. In this case we get
$X=\p^1 \times \nu_2(\p^1) \subset \p^1 \times \p^2 \subset \p^5$, 
which is not a \vs.

If $\delta_1=1$, we get $(2n_1-1)(2n_2-1)=9$ and,
if $n_1\leq n_2$,
$n_1=n_2=2$ (we then get $\p^2 \times \p^2$) or
$n_1=1$ and $n_2=5$; in this case the variety would be 
$\p^1 \times \p^5$ which is not a \vs.
\para

If $\cal G$ is simple,
we take a Borel subgroup ${\cal B} \subset {\cal G}$ and a maximal torus 
${\cal T} \subset {\cal B}$; let $\g$ be the Lie algebra associated to 
$\cal G$ and $\h$ the Cartan subalgebra associated to $\cal T$. Let then
$\Delta$ be the root system,
and $\Delta^+$ the set of all
positive roots, which by definition we ask to correspond to $\cal B$.
Let $\lambda$ be
the highest weight of $V$, 
$\mu=w_0(\lambda)$ the lowest weight, $l_{\lambda}$ and
$l_{\mu}$ the associated lines of eigenvectors. Let
$e_{\lambda}\in l_{\lambda}$ and $e_{\mu}\in l_{\mu}$ different from 0 and 
$s=e_{\lambda} + e_{\mu}$.
\begin{lemm}{$Sec(X)$ has finitely many $\cal G$-orbits, and
$s$ is in the open one.}
\end{lemm}

\dem By the Bruhat decomposition theorem, we see that $X \subset \p^m$ 
has only finitely many $\pi(\cal B)$-orbits so $Sec(X)$ has
only finitely many $\cal G$-orbits.
Let ${\cal B}_\lambda$ and ${\cal B}_\mu$ be the stabilizers of
 $l_{\lambda}$ and $l_{\mu}$. 
Let for any root $\alpha$, $\g_\alpha$ be the associated root space of $\g$.
Letting $\n^+=\sum_{\alpha \in \Delta^+} \g_\alpha$ and
$\n^-=\sum_{\alpha \in \Delta^+} \g_{-\alpha}$, we know that  
$\g = \h \oplus \n^+ \oplus \n^-$ and that the Lie algebras associated to
${\cal B}_\lambda$ and ${\cal B}_\mu$ are $\h \oplus \n^+$
and $\h \oplus \n^-$,
so that
${\cal B}_\lambda{\cal B}_\mu$ is a dense subset of $\cal G$. 
Thus ${\cal B}_\mu . l_{\lambda}={\cal B}_\mu{\cal B}_\lambda . l_{\lambda}$
contains an open subset $U$ of $X={\cal G}.l_\lambda$. Consequently, the
${\cal G}$-orbit of $l_\lambda \times l_\mu$ in $X \times X$ contains
$U \times l_\mu$ and so ${\cal G}.(U \times l_\mu)$, which is 
dense in $X \times X$.
\fin

Terracini's lemma (cf for example \cite{fulton}) allows us to conclude that\\
$T_sSec(X)=<T_{e_{\lambda}}X,T_{e_{\mu}}X>$.

As $X={\cal G}. l_{\lambda}$, denoting by $V_\alpha \subset V$ 
the weight subspace associated to the weight $\alpha$,
$$T_{e_{\lambda}}X=Ad(\g ) l_{\lambda}=l_\lambda \oplus 
(\oplus_{\alpha \in \Delta ^+} 
Ad(\g_{-\alpha}) l_{\lambda})
\subset 
l_\lambda \oplus
(\oplus_{\alpha \in \Delta ^+}
V_{\lambda - \alpha})$$
Similarly,
$$T_{e_\mu}X \subset l_\mu \oplus (\oplus_{\alpha \in \Delta ^+}
V_{\mu + \alpha})$$

As $\dim Sec(X) < 2\dim X +1$, there exists 
$\alpha,\beta \in \Delta^+$ such that $\lambda - \alpha=\mu + \beta$, that is
$\lambda - w_0(\lambda) = \alpha + \beta$. Let's notice that if $\lambda$
is the highest root, this equation has only one solution: 
$\alpha = \lambda$ and
$\beta = - w_0(\lambda) = \lambda$. 
The variety we get is then the closed orbit of the action of 
$\cal G$ in its projectiviced adjoint representation;
it is not a \vs as 
$T_{e_\lambda}X \cap T_{e_\mu}X$ equals the line 
$[\g_{\lambda},\g_{-\lambda}]$ of $V_0$ and so
$\dim Sec(X)=2 \dim X$.\\
In the general case, the existence of two positive roots 
whose sum is $\lambda - w_0(\lambda)$ is a very restrictive condition
which allows us to easily classify all homogeneous varieties whose
secant variety does not have maximal dimension. 
\para

I am now going to study two examples, a classical and an exceptional one,
of resolution of the equation 
$\lambda - w_0(\lambda) = \alpha + \beta$. I will use N. Bourbaki's notations
and results (cf \cite{bourbaki}).

\para
\begin{itemize}
  \item{{\bf root system of type $A_n$}\\
Here roots live in $\C^{n+1}$ and the positive roots
are the $\epsilon_i - \epsilon_j,1\leq i<j\leq n+1$. 
The fundamental weights are
$$\omega_i=(\epsilon_1 + \cdots + \epsilon_i) - \frac{i}{n+1}
\sum_{j=1}^{n+1} \epsilon_i$$ and $w_0=-Id$. Thus
$\omega_i - w_0(\omega_i) = \epsilon_1 + \cdots + \epsilon_i - 
\epsilon_{n+1-i} - \cdots - \epsilon_{n+1}$ if $i<\frac{n+1}{2}$
and $\omega_{n+1-i} - w_0(\omega_{n+1-i}) =\omega_i - w_0(\omega_i)$.
On the other hand, 
if $\alpha$ and $\beta$ are positive roots, and if
$\alpha + \beta = \sum s_i \epsilon_i$, then $\sum |s_i|\leq 4$.
As $\lambda$ is a sum of fundamental weights, we can deduce that only four
cases can occur:
\begin{itemize}
  \item{$\lambda = \omega_1$ or $\lambda = \omega_n$. Then 
$X=\p^n$.}
  \item{$\lambda = 2 \omega_1$ or $\lambda = 2 \omega_n$. Then
$X=\p^n \subset \p S^2\C^{n+1}$.}
  \item{$\lambda = \omega_2$ or $\lambda = \omega_{n-1}$. Then
$X=G(2,n+1)\subset \p \Lambda^2 \C^{n+1}$.}
  \item{$\lambda = \omega_1 + \omega_n$: adjoint representation.}
\end{itemize}}
\para
\item{{\bf root system of type $E_6$}\\
In this case, roots live in $\C^8$ and in the canonical basis 
$(\epsilon_i)$, the positive roots are the 
$\pm \epsilon_i + \epsilon_j \ (1\leq i<j\leq 5)$ and the
$\frac{1}{2}(\epsilon_8 - \epsilon_7 - \epsilon_6 + \sum_{i=1}^5 
(-1)^{\nu_i} \epsilon_i)$ with $\sum_{i=1}^5 \nu_i$ even.\\

Let's consider the application which sends $i$ to 
the sequence $(\lambda_j^i)_j$ such that if
$(\omega_i)$ are the fundamental weights, 
$\omega_i + w_0(\omega_i) = \sum_j \lambda_j^i \epsilon_j$. This application
is given by:

$$
\begin{array}{rcl}
1 & \mapsto & (0,0,0,0,1,-1,-1,1) \\
2 & \mapsto & (1,1,1,1,1,-1,-1,1) \\
3 & \mapsto & \frac{1}{2} (-1,1,1,3,3,-3,-3,3) \\
4 & \mapsto & (0,0,2,2,2,-2,-2,2) \\
5 & \mapsto & \frac{1}{2}(-1,1,1,3,3,-3,-3,3) \\
6 & \mapsto & (0,0,0,0,1,-1,-1,1)
\end{array}
$$

Given the formula for the positive roots,
$(\alpha + \beta)_8$ can only be $0,\frac{1}{2}$ or $1$ and as 
$\forall i, (\omega_i + w_0(\omega_i))_8>0 $, one can deduce 
that if $\lambda = \sum \lambda_i \omega_i$, only one 
$\lambda_i$ may be different from 0 and then equals 1. 
Finally, one has the list:
\begin{itemize}
  \item{$\lambda = \omega_1$ or $\lambda = \omega_6$: variety ${\cal E}_6$.}
  \item{$\lambda = \omega_2$: adjoint representation.}
\end{itemize}}

\end{itemize}
\para

One can treat all the cases this way, and check that the homogeneous
varieties whose secant variety does not have the maximal dimension are either
adjoint varieties, $X=\p^n$, quadrics, Veronese mappings
of $\p^n$ into $\p S^2\C^{n+1}$,
grassmannians of $2$-planes, possibly annihilating a quadratic or symplectic
form, or the closed orbit of the minimal representation 
of $E_6$,$F_4$ (a hyperplane section of that of $E_6$) or 
$G_2$ ($X$ is then also a quadric). 
An immediate computation of dimension then concludes the proof.
\fin

\section{More geometric properties}

We are now going to show that $G$ exchanges $X$ and $X^*$, and is a
birational  involution. I want to prove these properties without using
the previous classification.\\
As $D_{w_0}G$ is a linear isomorphism between $X$ and $Y$, $Y$ is a \vs.
So let $F^*$ be a trilinear symmetric form such that
$l\in Sec(Y) \Leftrightarrow F^*(l,l,l)=0$ and 
$l\in Y \Leftrightarrow F^*(l,l,.)=0$. Let
$G^*(l)=\frac{F^*(l,l,.)}{F^*(l,l,l)}$.

\begin{prop}{Let $w_0\in \C^{m+1}-Sec(X)$. Then $G^* \circ G (w_0)=w_0$.
\label{involution}}
\end{prop}

In particular $G$ defines a birational map from $\p^m$ 
to $\p^{m^*}$ with inverse
$G^*$.  The details of the proof
of this are left to the reader in \cite[p.79]{zak}
and written by
L. Ein and N. Shepherd-Barron (\cite[p.78,Theorem 2.6]{ein}).
They give two proofs but both of them require to know that $X$ is homogeneous;
and I pretend to give more elementary arguments.

\dem As $u\in Sec(X) \Leftrightarrow L_{w_0}(u)\in Sec(Y)$,
there exists $\lambda_{w_0}\in \C^*$ such that\\
$F^*[L_{w_0}(u),L_{w_0}(v),L_{w_0}(w)]=\lambda_{w_0}F(u,v,w)$  ($\diamond$)

Let $w_0\in \C^{m+1}-Sec(X)$ and $l\in {(\C^{m+1})}^*$, we want to see that
$$\frac{l(w_0)}{F(w_0)}=\frac{F^*(w_0^*,w_0^*,l)}{F^*(w_0^*)}$$
As $L_{w_0}(w_0)=-F(w_0)w_0^*$, applying ($\diamond$) with
$u=v=w_0$ and $w$ such that $L_{w_0}(w)=l$ yields:\\
$F^*(w_0^*,w_0^*,l) F^2(w_0)=\lambda_{w_0}F(w_0,w_0,w)$. If we take 
$u=v=w=w_0$, we get:\\
$F^*(w_0^*,w_0^*,w_0^*)F^3(w_0)=-\lambda_{w_0}F(w_0)$.
The quotient of the last two equalities yields the expected equality.
\fin

$G$ considered as a rational morphism can be extended to 
$\se$ by $G(p)=p^*$ and then $G(\se)=Y$.
Let's denote by $G(X)$ the total transform of $X$, that is the set of 
limits of $G(x_n)$ for
$(x_n)$ any sequence converging to an element of $X$. 
The next proposition is stated
without proof in \cite[p.79]{zak}:

\begin{prop}{$G(X)=Sec(Y)=X^*$.}
\end{prop}

\dem We already know that $G(X)\subset Sec(Y)$ since out of $Sec(Y)$, $G$ is
invertible. In the other direction, if $y_n\rightarrow y$ with $y\in Sec(Y)-Y$
and $y_n \not \in Sec(Y)$, then letting $x_n=G^*(y_n)$,
proposition \ref{involution} yields $G(x_n)=y_n$, so $G(x_n)\rightarrow y$. As
$x_n\rightarrow G(y)\in X$, $y\in G(X)$, and $G(X)=Sec(Y)$.

As far as the second equality is concerned, let $w\in \p^m-Sec(X)$ and 
$w'=G(w)$. We have isomorphisms $D_wG$ and $D_{w'}G^*$ respectively between
$X$ and $Y$, and between $Y$ and ${(Sec(Y))}^*$. So we get a composed 
isomorphism
$X \simeq {(Sec(Y))}^*$. But by the choice of $w$ and $w'$ and
proposition \ref{involution}, this composed isomorphism is the identity,
so that
$X={(Sec(Y))}^*$ and $X^*=Sec(Y)$.  \fin
\para

We now show that $\se$ is homogeneous.
This implies that under the action of $\cal G$,
$\C^{m+1}$ has only one invariant, $F$, meaning that every $\cal G$-invariant
polynomial is, up to a constant, a power of $F$. 
We can then use the classification of 
all groups whose invariant algebra is free (cf \cite{kac}) to deduce another
proof of the classification of Severi varieties (theorem
\ref{classification}).

Let $p\in \se$ and
$P(w_0)=2F(w_0)p-6w_0^*(p)w_0$.
We start with two technical lemmas:

\begin{lemm}{There exists $w_0$ such that $P(w_0)\not\in Sec(X)$.}
\end{lemm}

\dem If $F[P(w)]=0$ were true for all $w\in \C^{m+1}$, then we would get, 
either $F(w)^2F(w,w,p)F(w,p,p)=0$ for all $w$ and so
$F(w,w,p)=0$ for all $w$, either $F(w,p,p)=0$ for all $w$. In the first
case, polarisation yields also $F(w,p,p)=0$; so in both cases 
$p\in X$, contradicting the choice of $p$. \fin

\begin{lemm}{Let $\omega\not \in Sec(X)$. The matrix $F(\omega,.,.)$ is
inversible.  \label{matinversible}}
\end{lemm}

\dem This matrix is the differential in $\omega$ 
of $\fonc{H}{w}{F(w,w,.)}$. As $F(\omega)\not =0$, 
we can define in a neigbourhood of
$\omega$ a square root $\sqrt{F(w)}$, wich is not 0, and proposition
\ref{involution} yields that the composition
$$w\mapsto w_0:=\frac{w}{\sqrt{F(w)}} \stackrel{H}{\mapsto}
F(w_0,w_0,.) \stackrel{G^*}{\mapsto} G^*[F(w_0,w_0,.)]$$
is the identity, so that $F(\omega,.,.)$ is inversible.
\fin

\begin{prop}{$\se$ is homogeneous.  \label{homo} }
\end{prop}

\dem Let $p\in \se$ and $w_0$ such that $P(w_0)\not \in Sec(X)$.
Let $L(w)$ the derivative of $L_{\omega}(p)$ when $\omega$ tends
to $w_0$ in the direction of $w$:
$$
\begin{array}{ll}
L(w) & =6F(w_0,w_0,w)F(w_0,p,.)+2F(w_0,w_0,w_0)F(w,p,.)\\
     & -6F(w_0,w,p)F(w_0,w_0,.)-6F(w_0,w_0,p)F(w_0,w,.)
\end{array}
$$
I'm going to explain why Ker$(L)\cap$ Ker$(w_0^*)= \{ 0 \}$, 
so that $L$ has rank at least $m$. Therefore, if $p$ and $p'$
are elements of $\se$,
$\{L_w(p)\}$ and $\{L_w(p')\}$ contain open subsets of $Sec(Y)$ so their
intersection is not empty, concluding the proof of the proposition.

So let $w$ such that $L(w)=0$ and $F(w_0,w_0,w)=0$.  $L(w).w_0=0$
yields $F(w_0,w,p)=0$; as $F(w_0,w_0,w)=0$, we deduce
$F[P(w_0),w,.]=0$. As by definition of $w_0$, $P(w_0)\not \in Sec(X)$,
lemma \ref{matinversible} yields $w=0$. \fin

\vfill
\begin{flushleft}
Pierre-Emmanuel CHAPUT\\
INSTITUT FOURIER\\
Laboratoire de Math\'ematiques\\
UMR5582 (UJF-CNRS)\\
BP 74\\
38402 St MARTIN D'HÈRES Cedex (France)\\
\vspace{.5cm} chaput@ujf-grenoble.fr
\end{flushleft}

\end{document}